# MULTISCALE INFERENCE ABOUT A DENSITY

By Lutz Dümbgen[1] and Günther Walther[2]

*University of Bern and Stanford University*

We introduce a multiscale test statistic based on local order statistics and spacings that provides simultaneous confidence statements for the existence and location of local increases and decreases of a density or a failure rate. The procedure provides guaranteed finite-sample significance levels, is easy to implement and possesses certain asymptotic optimality and adaptivity properties.

**1. Introduction.** An important aspect in the analysis of univariate data is inference about qualitative characteristics of their distribution function $F$ or density $f$, such as the number and location of monotone or convex regions, local extrema or inflection points. This issue has been addressed in the literature using a variety of methods. Silverman (1981), Mammen, Marron and Fisher (1992), Minnotte and Scott (1993), Fisher, Mammen and Marron (1994), Minnotte (1997), Cheng and Hall (1999) and Chaudhuri and Marron (1999, 2000) use kernel density estimates. Excess masses and related ideas are employed by Hartigan and Hartigan (1985), Hartigan (1987), Müller and Sawitzky (1991), Polonik (1995) and Cheng and Hall (1998). Good and Gaskins (1980) and Walther (2001) use maximum likelihood methods, whereas Davies and Kovac (2004) employ the taut string method. In the present paper, a qualitative analysis of a density $f$ means simultaneous confidence statements about regions of increase and decrease, as well as local extrema. Such simultaneous inference has been only sparingly treated in the literature. Also, the methods available thus far provide only approximate significance levels as the sample size tends to infinity, and they rely on certain regularity conditions on $f$.

Received February 2006; revised June 2007.
[1]Supported by Swiss National Science Foundation.
[2]Supported by NSF Grants DMS-98-75598, DMS-05-05682 and NIH Grant 5R33HL068522.
*AMS 2000 subject classifications.* 62G07, 62G10, 62G15, 62G20, 62G30.
*Key words and phrases.* Exponential inequality, modes, monotone failure rate, multiple test, order statistics, spacings, subexponential increments, sub-Gaussian tails.

In this paper, we introduce and analyze a procedure that provides simultaneous confidence statements with guaranteed given significance level for arbitrary sample size. The approach is similar to that of Dümbgen (2002), who used local rank tests in the context of nonparametric regression, or Chaudhuri and Marron's (1999, 2000) SiZer, where kernel estimators with a broad range of bandwidths are combined. Here, we utilize test statistics based on local order statistics and spacings. The use of spacings for nonparametric inference about densities has a long history. For instance, Pyke (1995) describes various goodness-of-fit tests based on spacings and Roeder (1992) uses such tests for inference about normal mixtures. Confidence bands for an antitonic density on $[0, \infty)$ via uniform order statistics and spacings have been constructed by Hengartner and Stark (1995) and Dümbgen (1998).

In Section 2, we define local spacings and related test statistics which indicate isotonic or antitonic trends of $f$ on certain intervals. Then, a deterministic inequality (Proposition 2.1) relates the joint distribution of all these test statistics in general to the distribution in the special case of a uniform density. This enables us to define a multiple test concerning monotonicity properties of $f$. Roughly speaking, we consider all intervals whose endpoints are observations. The rationale for using and combining statistics corresponding to such a large collection of (random) intervals is that the power for detecting an increase or decrease of $f$ is maximized when the tested interval is close to an interval on which $f$ has such a trend. In that context, we also discuss two important differences with Chaudhuri and Marron's SiZer map.

In Section 3, we describe a particular way of calibrating and combining the single test statistics. Optimality results in Section 4 show that in many relevant situations, the resulting multiscale test is asymptotically as powerful, in the minimax sense, as any procedure can essentially be for detecting increases and decreases of $f$ on small intervals as well as on large intervals. Thus, neither the guaranteed confidence level nor the simultaneous consideration of many intervals results in a substantial loss of power. In addition, we prove that our procedure is able to detect and localize an arbitrary number of local extrema under weak assumptions on the strength of these effects.

In Section 5, we consider a density $f$ on $(0, \infty)$ and modify our multiple test in order to analyze monotonicity properties of the failure rate $f/(1 - F)$. It is well known that spacings are a useful object in this context; see, for example, Proschan and Pyke (1967), Bickel and Doksum (1969) and Barlow and Doksum (1972). While these authors use global test statistics, Gijbels and Heckman (2004) localize, standardize and combine such tests, albeit without calibrating the various scales. Hall and van Keilegom (2005) use resampling from an appropriately calibrated null distribution in order to achieve better sensitivity in detecting local effects, which leads to an asymptotically valid test procedure without explicit information about the location of these effects. Walther (2001) uses a multiscale maximum likelihood analysis to detect local effects.

Section 6 illustrates the multiscale procedures with two examples and introduces a graphical display. Proofs and technical arguments are deferred to Section 7. One essential ingredient is an auxiliary result concerning stochastic processes with sub-Gaussian marginals and subexponential increments. This result generalizes Theorem 6.1 of Dümbgen and Spokoiny (2001) and is a corollary to more general results in a technical report by Dümbgen and Walther (2006).

To establish notation for the sequel, suppose that $Y_1, Y_2, \ldots, Y_m$ are independent random variables with unknown distribution function $F$ and (Lebesgue) density $f$ on the real line. In order to infer properties of $f$ from these data, we consider the corresponding order statistics $Y_{(1)} < Y_{(2)} < \cdots < Y_{(m)}$. In some applications, $F$ is known to be supported by an interval $[a, \infty)$, $(-\infty, b]$ or $[a, b]$, where $-\infty < a < b < \infty$. In that case, we add the point $Y_{(0)} := a$, the point $Y_{(m+1)} := b$ or both respectively, to our ordered sample. This yields a data vector $\mathbf{X} = (X_{(i)})_{i=0}^{n+1}$ with real components $X_{(0)} < X_{(1)} < \cdots < X_{(n+1)}$, where $n \in \{m-2, m-1, m\}$. For $0 \le j < k \le n+1$ with $k - j > 1$, the conditional joint distribution of $X_{(j+1)}, \ldots, X_{(k-1)}$, given $X_{(j)}$ and $X_{(k)}$, coincides with the joint distribution of the order statistics of $k - j - 1$ independent random variables with density

$$f_{jk}(x) := \frac{1\{x \in \mathcal{I}_{jk}\}f(x)}{F(X_{(k)}) - F(X_{(j)})},$$

where $\mathcal{I}_{jk}$ stands for the interval

$$\mathcal{I}_{jk} := (X_{(j)}, X_{(k)}).$$

Thus, $(X_{(j+i)})_{i=0}^{k-j}$ is useful for inferring properties of $f$ on $\mathcal{I}_{jk}$. The multiple tests which follow are based on all such tuples.

**2. Local spacings and monotonicity properties of the density.** Let us consider one particular interval $\mathcal{I}_{jk}$ and condition on its endpoints. In order to test whether $f$ is nonincreasing or nondecreasing on $\mathcal{I}_{jk}$, we introduce the local order statistics

$$X_{(i;j,k)} := \frac{X_{(i)} - X_{(j)}}{X_{(k)} - X_{(j)}}, \qquad j \le i \le k,$$

and the test statistic

$$T_{jk}(\mathbf{X}) := \sum_{i=j+1}^{k-1} \beta(X_{(i;j,k)}),$$

where

$$\beta(x) := 1\{x \in (0,1)\}(2x - 1).$$

This particular test statistic $T_{jk}(\mathbf{X})$ appears as a locally most powerful test statistic for the null hypothesis "$\lambda \le 0$" versus "$\lambda > 0$" in the parametric model, where

$$f_{jk}(x) = \frac{1\{x \in \mathcal{I}_{jk}\}}{X_{(k)} - X_{(j)}}\left(1 + \lambda\left(\frac{x - X_{(j)}}{X_{(k)} - X_{(j)}} - \frac{1}{2}\right)\right).$$

Elementary algebra yields the following alternative representation of our single test statistics:

$$(2.1) \quad T_{jk}(\mathbf{X}) = -(k-j)\sum_{i=j+1}^{k} \beta\left(\frac{i - j - 1/2}{k - j}\right)(X_{(i;j,k)} - X_{(i-1;j,k)}).$$

Thus, $T_{jk}(\mathbf{X})$ is a weighted average of the local spacings $X_{(i;j,k)} - X_{(i-1;j,k)}$, $j < i \le k$.

Suppose that $f$ is constant on $\mathcal{I}_{jk}$. The random variable $T_{jk}(\mathbf{X})$ is then distributed (conditionally) as

$$(2.2) \quad \sum_{i=1}^{k-j-1} \beta(U_i),$$

with independent random variables $U_i$ having uniform distribution on $[0,1]$. Note that the latter random variable has mean zero and variance $(k - j - 1)/3$. However, if $f$ is nondecreasing or nonincreasing on $\mathcal{I}_{jk}$, then $T_{jk}(\mathbf{X})$ tends to be positive or negative, respectively. The following proposition provides a more general statement, which is the key to our multiple test.

PROPOSITION 2.1. *Define* $\mathbf{U} = (U_{(i)})_{i=0}^{n+1}$, *with components* $U_{(i)} := F_o(X_{(i)})$, *where* $F_o$ *is the distribution function corresponding to the density* $f_{0,n+1}$. *Then,* $U_{(1)}, \ldots, U_{(n)}$ *are distributed as the order statistics of* $n$ *independent random variables having uniform distribution on* $[0,1]$, *while* $U_{(0)} = 0$ *and* $U_{(n+1)} = 1$. *Moreover, for arbitrary integers* $0 \le j < k \le n+1$ *with* $k - j > 1$,

$$T_{jk}(\mathbf{X}) \begin{cases} \ge T_{jk}(\mathbf{U}), & \text{if } f \text{ is nondecreasing on } \mathcal{I}_{jk}, \\ \le T_{jk}(\mathbf{U}), & \text{if } f \text{ is nonincreasing on } \mathcal{I}_{jk}. \end{cases}$$

This proposition suggests the following multiple test. Suppose that for a given level $\alpha \in (0,1)$, we know constants $c_{jk}(\alpha)$ such that

$$(2.3) \quad \mathbb{P}\{|T_{jk}(\mathbf{U})| \le c_{jk}(\alpha) \text{ for all } 0 \le j < k \le n+1, k - j > 1\} \ge 1 - \alpha.$$

Let

$$\mathcal{D}^{\pm}(\alpha) := \{\mathcal{I}_{jk} : \pm T_{jk}(\mathbf{X}) > c_{jk}(\alpha)\}.$$

One can then claim with confidence $1 - \alpha$ that $f$ must have an increase on every interval in $\mathcal{D}^+(\alpha)$ and that it must have a decrease on every interval in $\mathcal{D}^-(\alpha)$. In other words, with confidence $1 - \alpha$, we may claim that for every $\mathcal{I} \in \mathcal{D}^\pm(\alpha)$ and for every version of $f$, there exist points $x, y \in \mathcal{I}$ with $x < y$ and $\pm(f(y) - f(x)) > 0$.

Combining the two families $\mathcal{D}^\pm(\alpha)$ properly allows the detection and localization of extrema as well. Suppose, for instance, that there exist intervals $I_1, I_2, \ldots, I_m$ in $\mathcal{D}^+(\alpha)$ and $D_1, D_2, \ldots, D_m$ in $\mathcal{D}^-(\alpha)$ such that $I_1 \leq D_1 \leq I_2 \leq D_2 \leq \cdots \leq I_m \leq D_m$, where the inequalities are to be understood elementwise. Under the weak assumption that $f$ is continuous, one can conclude with confidence $1 - \alpha$ that $f$ has at least $m$ different local maxima and $m - 1$ different local minima.

Note that our multiscale test allows the combination of test statistics $T_{jk}(\mathbf{X})$ with arbitrary "scales" $k - j$. This is an advantage over Chaudhuri and Marron's (1999, 2000) SiZer map, where statements about *multiple* increases and decreases are available only at a common bandwidth. This is due to the fact that these authors use kernels with unbounded support and rely on a particular variation-reducing property of the Gaussian kernel which holds only for an arbitrary but global bandwidth. Another consequence of the kernel's unbounded support is that localizing trends of $f$ itself is not possible.

**3. Properly combining the single test statistics.** It remains to define constants $c_{jk}(\alpha)$ satisfying (2.3). First, note that $T_{jk}(\mathbf{U})$ has mean zero and standard deviation $\sqrt{(k-j-1)/3}$. Motivated by recent results of Dümbgen and Spokoiny (2001) concerning multiscale testing in Gaussian white noise models, we consider the test statistic

$$T_n(\mathbf{X}) := \max_{0 \leq j < k \leq n+1 : k-j > 1} \left( \sqrt{\frac{3}{k-j-1}} |T_{jk}(\mathbf{X})| - \Gamma\left(\frac{k-j}{n+1}\right) \right),$$

where $\Gamma(\delta) := (2 \log(e/\delta))^{1/2}$. This particular additive calibration for various scales is necessary for the optimality results to follow. Without the term $\Gamma((k-j)/(n+1))$, the null distribution would be dominated by small scales, as there are many more local test statistics on small scales than on large scales, with a corresponding loss of power at large scales. The next theorem states that our particular test statistic $T_n(\mathbf{U})$ converges in distribution. Unless stated otherwise, asymptotic statements in this paper refer to $n \to \infty$.

THEOREM 3.1.

$$T_n(\mathbf{U}) \to_{\mathcal{L}} T(W) := \sup_{0 \leq u < v \leq 1} \left( \frac{|Z(u,v)|}{\sqrt{v-u}} - \Gamma(v-u) \right),$$

*where*

$$Z(u,v) := 3^{1/2} \int_u^v \beta\left(\frac{x-u}{v-u}\right) dW(x)$$

*and $W$ is a standard Brownian motion on $[0,1]$. Moreover, $0 \le T < \infty$ almost surely.*

Consequently, if $\kappa_n(\alpha)$ denotes the $(1-\alpha)$-quantile of $\mathcal{L}(T_n(\mathbf{U}))$, then $\kappa_n(\alpha) = O(1)$ and the constants

$$c_{jk}(\alpha) := \sqrt{\frac{k-j-1}{3}} \left( \Gamma\left(\frac{k-j}{n+1}\right) + \kappa_n(\alpha) \right)$$

satisfy requirement (2.3). For explicit applications, we do not use the limiting distribution in Theorem 3.1, but rely on Monte Carlo simulations of $T_n(\mathbf{U})$ which are easily implemented.

**4. Power considerations.** Throughout this section, we focus on the detection of increases of $f$ by means of $\mathcal{D}^+(\alpha)$. Analogous results hold true for decreases of $f$ and $\mathcal{D}^-(\alpha)$.

For any bounded open interval $I \subset \mathbb{R}$, we quantify the isotonicity of $f$ on $I$ by

$$\inf_I f' := \inf_{x,y \in I: x<y} \frac{f(y)-f(x)}{y-x}$$
$$= \inf_{x \in I} f'(x) \quad \text{if } f \text{ is differentiable on } I.$$

We now analyze the difficulty of detecting intervals $I$ with $\inf_I f' > 0$. An appropriate measure of this difficulty turns out to be

$$H(f,I) := \inf_I f' \cdot |I|^2 / \sqrt{F(I)},$$

where $|I|$ denotes the length of $I$. Note that this quantity is affine equivariant, in the sense that it does not change when $f$ and $I$ are replaced by $\sigma^{-1} f(\sigma^{-1}(\cdot - \mu))$ and $\{\mu + \sigma x : x \in I\}$, respectively, with $\mu \in \mathbb{R}$, $\sigma > 0$. For given numbers $\delta \in (0,1]$ and $\eta \in \mathbb{R}$, we define

$$\mathcal{F}(I,\delta,\eta) := \{f : F(I) = \delta, H(f,I) \ge \eta\}$$

and

$$\mathcal{F}(\delta,\eta) := \bigcup_{\text{bounded intervals } I} \mathcal{F}(I,\delta,\eta).$$

Note that $f(x) \ge \inf_I f' \cdot (x - \inf(I))$ on $I$, so $F(I) \ge \inf_I f' \cdot |I|^2/2$. Hence,

(4.1) $$H(f,I) \le 2\sqrt{F(I)}.$$

Thus, $\mathcal{F}(I,\delta,\eta)$ and $\mathcal{F}(\delta,\eta)$ are nonvoid if and only if $\eta \le 2\sqrt{\delta}$.

THEOREM 4.1. *Let $\delta_n \in (0,1]$ and $0 < c_n < \sqrt{24} < C_n$.*

(a) *Let $I_n$ be a bounded interval and $f_n$ a density in $\mathcal{F}(I_n, \delta_n, C_n \times \sqrt{\log(e/\delta_n)/n})$. Then,*

$$\mathbb{P}_{f_n}(\mathcal{D}^+(\alpha) \text{ contains an interval } J \subset I_n) \to 1,$$

*provided that $(C_n - \sqrt{24})\sqrt{\log(e/\delta_n)} \to \infty$.*

(b) *Let $\phi_n(\mathbf{X})$ be any test with level $\alpha \in (0,1)$ under the null hypothesis that $\mathbf{X}$ is drawn from a nonincreasing density. If $(\log n)^2/n \le \delta_n \to 0$, then*

$$\inf_{f \in \mathcal{F}(\delta_n, c_n\sqrt{\log(e/\delta_n)/n})} \mathbb{E}_f \phi_n(\mathbf{X}) \le \alpha + o(1),$$

*provided that $(\sqrt{24} - c_n)\sqrt{\log(e/\delta_n)} \to \infty$.*

(c) *Let $I_n$ be any interval and $b_n$ some number in $[0, 2\sqrt{n\delta_n}]$. If $\phi_n(\mathbf{X})$ is any test with level $\alpha \in (0,1)$ under the null hypothesis that the density is nonincreasing on $I_n$, then*

$$\inf_{f \in \mathcal{F}(I_n, \delta_n, b_n/\sqrt{n})} \mathbb{E}_f \phi_n(\mathbf{X}) \to 1$$

*implies that $b_n \to \infty$ and $n\delta_n \to \infty$.*

Theorem 4.1 establishes that our multiscale statistic is optimal, in the asymptotic minimax sense, for detecting an increase on an unknown interval, both in the case of an increase occurring on a small scale ($\delta_n \searrow 0$) and in the case of an increase occurring on a large scale ($\liminf \delta_n > 0$).

In the case of small scales, a comparison of (a) and (b) shows that there is a cut-off for the quantity $H(f, I)$ at $\sqrt{24 \log(e/\delta_n)/n}$: if one replaces the factor 24 with $24 + \epsilon_n$, with $\epsilon_n \searrow 0$ sufficiently slowly, then the multiscale test will detect and localize such an increase with asymptotic power one, whereas in the case $24 - \epsilon_n$, no procedure can detect such an increase with nontrivial asymptotic power.

In the case of large scales, one may replace $\mathcal{F}(I_n, \delta_n, C_n\sqrt{\log(e/\delta_n)/n})$ in (a) with the family $\mathcal{F}(I_n, \delta_n, \tilde{C}_n/\sqrt{n})$, where $\tilde{C}_n \to \infty$. Then, a comparison of (a) and (c) again shows our multiscale test to be optimal, even in comparison to tests using a priori knowledge of the location and scale of the potential increase. Hence, searching over all (large and small) scales does not involve a serious drawback. In the case of small scales, (a) and (c) together show that ignoring prior information about the location of the potential increase leads to a penalty factor of order $o(\sqrt{\log(e/\delta_n)}) = o(\sqrt{\log n})$.

EXAMPLE 1. Let us first illustrate the theorem in the special case of a fixed continuous density $f$ and a sequence of intervals $I_n$ converging to a given point $x_o$, where we use the abbreviation

$$\rho_n := \log(n)/n.$$

EXAMPLE 1A. Let $f$ be continuously differentiable in a neighborhood of $x_o$, such that $f(x_o) > 0$ and $f'(x_o) > 0$. If $|I_n| = D_n \rho_n^{1/3}$ with $D_n \to D > 0$, then $\delta_n := F(I_n)$ is equal to $D_n f(x_o) \rho_n^{1/3}(1 + o(1))$ and $\inf_{I_n} f' = f'(x_o) + o(1)$. Hence, the quantity $H(f, I_n)$ may be written as $D_n^{3/2} f'(x_o) f(x_o)^{-1/2} \times \rho_n^{1/2}(1 + o(1))$, while $\sqrt{24 \log(e/\delta_n)/n} = 8^{1/2} \rho_n^{1/2} + o(1)$. Consequently, the conclusion of Theorem 4.1(a) is correct if

$$D_n \searrow (8f(x_o)/f'(x_o)^2)^{1/3}$$

sufficiently slowly.

EXAMPLE 1B. Let $f$ be differentiable on $(x_o, \infty)$, with $f(x_o) = 0$ and $f'(x_o + h) = \gamma h^{\kappa - 1}(1 + o(1))$ as $h \searrow 0$, where $\gamma, \kappa > 0$. Defining the interval $I_n$ to be $[x_o + C_1 \rho_n^{1/(\kappa+1)}, x_o + C_2 \rho_n^{1/(\kappa+1)}]$ with $0 \le C_1 < C_2$, the conclusion of Theorem 4.1(a) is correct, provided that $\min(C_1^{\kappa-1}, C_2^{\kappa-1})$ and $C_2/C_1$ are sufficiently large.

EXAMPLE 1C. Let $f$ be twice continuously differentiable in a neighborhood of $x_o$, such that $f(x_o) > 0$, $f'(x_o) = 0$ and $\pm f''(x_o) \ne 0$. Now, take the two intervals $I_n^{(\ell)} := [x_o - C_2 \rho_n^{1/5}, x_o - C_1 \rho_n^{1/5}]$ and $I_n^{(r)} := [x_o + C_1 \rho_n^{1/5}, x_o + C_2 \rho_n^{1/5}]$, with $0 < C_1 < C_2$. If $C_1$ and $C_2/C_1$ are sufficiently large, then it follows from Theorem 4.1(a) and its extension to locally decreasing densities that

$$\mathbb{P}(\mathcal{D}^\pm \text{ contains some } J \subset I_n^{(\ell)} \text{ and } \mathcal{D}^\mp \text{ contains some } J \subset I_n^{(r)}) \to 1.$$

Thus, our multiscale procedure will detect the presence of the mode with asymptotic probability one and furthermore localize it with precision $O_p((\log(n)/n)^{1/5})$. Up to the logarithmic factor, this is the optimal rate for estimating the mode [cf. Hasminskii (1979)].

EXAMPLE 2. Now, let $I$ be a fixed bounded interval and consider a sequence of densities $f_n$ such that $\sup_{x \in I} |f_n(x) - f_o| \to 0$ for some constant $f_o > 0$. Here, the conclusion of Theorem 4.1(a) is correct, provided that

$$\sqrt{n} \cdot \inf_I f_n' \to \infty.$$

The next theorem concerns the simultaneous detection of several increases of $f$.

THEOREM 4.2. *Let* $f = f_n$ *and let* $\mathcal{I}_n$ *be a collection of nonoverlapping bounded intervals such that for each* $I \in \mathcal{I}_n$,

(4.2) $$H(f_n, I) \ge C(\sqrt{\log(e/F_n(I))} + b_n)/\sqrt{n},$$

with constants $0 \le b_n \to \infty$ and $C \ge \sqrt{24}$. Then,

$$\mathbb{P}_{f_n}(\text{for each } I \in \mathcal{I}_n, \mathcal{D}^+(\alpha) \text{ contains an interval } J \subset I) \to 1$$

in each of the following three settings, where $\delta_n := \min_{I \in \mathcal{I}_n} F_n(I)$:

(i) $C \ge 34$;
(ii) $C > 2\sqrt{24}$ and $n\delta_n/\log(e\#\mathcal{I}_n) \to \infty$;
(iii) $C = \sqrt{24}$ and $n\delta_n/\log(e\#\mathcal{I}_n) \to \infty$, $\log \#\mathcal{I}_n = o(b_n^2)$.

It will be shown in Section 7 that (4.2) entails $n\delta_n \ge (C^2/4 + o(1))\log n$. In particular, $\#\mathcal{I}_n = o(n)$. Moreover, Theorem 4.1(a) follows from Theorem 4.2 by considering setting (iii) with $\mathcal{I}_n$ consisting of a single interval $I_n$.

A comparison with Theorem 4.1(a) shows that the price for the simultaneous detection of an increasing number of increases or decreases is essentially a potential increase of the constant $\sqrt{24}$.

The proof of Theorem 4.2 rests on an inequality involving the following auxiliary functions. For $c \in [-2, 2]$ and $u \in [0, 1]$, let

$$g_c(u) := 1 + c(u - 1/2).$$

This defines a probability density on $[0, 1]$ with distribution function

$$G_c(u) := u - cu(1-u)/2.$$

PROPOSITION 4.1. *Define* $\mathbf{U} = (U_{(i)})_{i=0}^{n+1}$ *as in Proposition 2.1. For arbitrary integers* $0 \le j < k \le n+1$ *with* $k - j > 1$, *it follows from* $\inf_{\mathcal{I}_{jk}} f' \ge 0$ *that*

$$T_{jk}(\mathbf{X}) \ge \sum_{i=j+1}^{k-1} \beta(G_S^{-1}(U_{(i;j,k)})) \qquad \text{with } S := \frac{H(f, \mathcal{I}_{jk})}{\sqrt{F(\mathcal{I}_{jk})}}.$$

*Moreover, for any fixed* $c \in [-2, 2]$ *and* $U \sim \text{Unif}[0, 1]$,

$$\mathbb{E}\beta(G_c^{-1}(U)) = c/6, \qquad \text{Var}(\beta(G_c^{-1}(U))) \le 1/3,$$

*while*

$$\mathbb{E}\exp(t\beta(G_c^{-1}(U))) \le \exp(ct/6 + t^2/6) \qquad \text{for all } t \in \mathbb{R}.$$

**5. Monotonicity of the failure rate.** To investigate local monotonicity properties of the failure rate $f/(1-F)$, such as the presence of a 'burn-in' period or a 'wear-out' period, we consider

$$W_i := \sum_{k=1}^{i} D_k \Big/ \sum_{k=1}^{n+1} D_k, \qquad i = 0, \ldots, n+1,$$

where $D_i := (n - i + 2)(X_{(i)} - X_{(i-1)})$, $i = 1, \ldots, n+1$, are the normalized spacings. Here, $X_{(0)} < X_{(1)} < \cdots < X_{(n+1)}$ are the order statistics of $n + 2$ or $n + 1$ i.i.d. observations from $F$, in the latter case with $X_{(0)}$ being the left endpoint of the support of $F$. The next proposition shows that the problem can now be addressed by applying the methodology of Section 2 to the transformed data vector $\mathbf{W} = (W_i)_{i=0}^{n+1}$.

PROPOSITION 5.1. *Let* $X'_{(i)} := -\log(1 - F(X_{(i)}))$, $i = 0, \ldots, n+1$, *and define* $\mathbf{W}' = (W'_i)_{i=0}^{n+1}$ *analogously to above, with* $\mathbf{X}'$ *in place of* $\mathbf{X}$. *Then,* $\mathbf{W}' =_{\mathcal{L}} \mathbf{U}$ *and, for arbitrary integers* $0 \le j < k \le n+1$ *with* $k - j > 1$,

$$T_{jk}(\mathbf{W}) \begin{cases} \ge T_{jk}(\mathbf{W}'), & \textit{if the failure rate of } f \textit{ is nondecreasing on } \mathcal{I}_{jk}, \\ \le T_{jk}(\mathbf{W}'), & \textit{if the failure rate of } f \textit{ is nonincreasing on } \mathcal{I}_{jk}. \end{cases}$$

**6. Graphical displays and examples.** We first illustrate our methodology with a sample of size $m = 300$ from the mixture distribution

$$F = 0.3 \cdot \mathrm{Gamma}(2) + 0.2 \cdot \mathcal{N}(5, 0.1) + 0.5 \cdot \mathcal{N}(11, 9),$$

where Gamma(2) denotes the gamma distribution with density $g(x) = xe^{-x}$ on $(0, \infty)$. Figure 1 depicts the density $f$ of $F$.

Figure 2 provides a line plot of the data and a visual display of the multiscale analysis. The horizontal line segments above the line plot depict all minimal intervals in $\mathcal{D}^+(0.1)$, while those below the line plot depict all minimal intervals in $\mathcal{D}^-(0.1)$. Here, we estimated the quantile $\kappa_{m-2}(0.1)$ to be 1.518 in 9999 Monte Carlo simulations, where we restricted $(j,k)$ in the definition of $T$ to index pairs $(j,k)$ such that $(k-j)/(m+1) \leq 0.34$. For example, we can conclude with simultaneous confidence 90% that each of the intervals $(0.506, 3.887)$ and $(5.022, 5.841)$ contains a decrease and that each of the intervals $(3.983, 4.882)$ and $(5.841, 10.307)$ contains an increase. As these four intervals are disjoint, we can conclude with confidence 90% that the density has at least three modes.

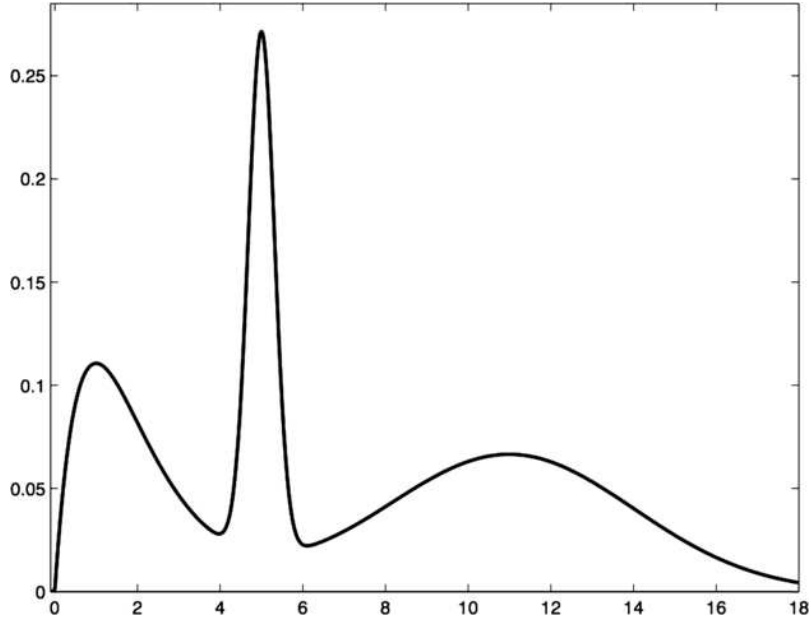

FIG. 1. *Density of* $0.3 \cdot \mathrm{Gamma}(2) + 0.2 \cdot \mathcal{N}(5, 0.1) + 0.5 \cdot \mathcal{N}(11, 9)$.

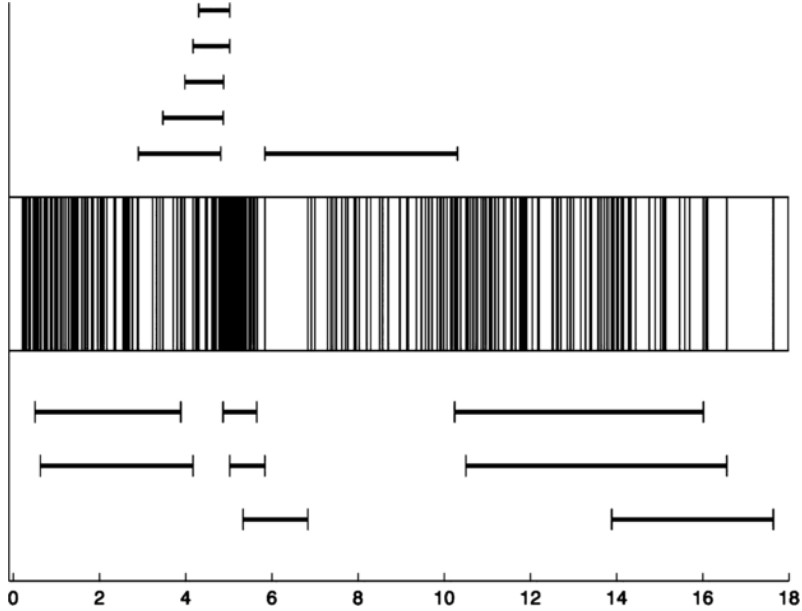

FIG. 2. *Minimal intervals in $\mathcal{D}^+(0.1)$ (top) and $\mathcal{D}^-(0.1)$ (bottom).*

A referee reports that the taut string method of Davies and Kovac (2004) found three modes in about 82% of cases. Our method finds three modes in about 39% and exactly two modes in about 50% of the cases. However, the latter method also allows the localization of the modes. Figure 3 provides a diagnostic tool for this type of inference. Each horizontal line segment, labeled by "+" or "−", depicts an interval in some $\mathcal{D}^+(\alpha)$ [resp., $\mathcal{D}^-(\alpha)$]. In each row, the depicted intervals are disjoint, with an alternating sequence of signs. The number in the first column gives the smallest significance level at which this sequence of alternating signs occurs, and the plot shows all such sequences that have a significance level of 10% or less. The intervals depicted in a given row are chosen to have the smallest right endpoint among the minimal intervals at the stated level. Consecutive intervals are plotted with a small vertical offset to more readily show their endpoints. For example,

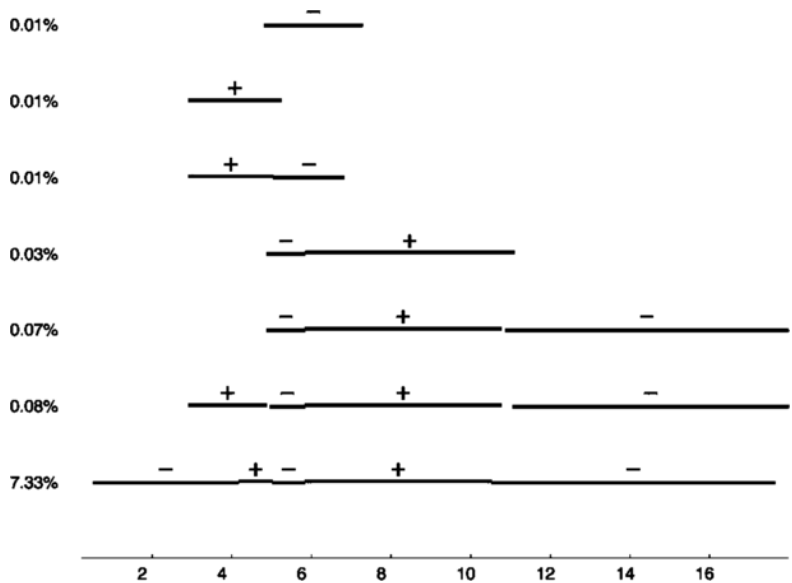

Fig. 3. *Alternating sequences of minimal intervals in $\mathcal{D}^+(\alpha)$, and $\mathcal{D}^-(\alpha)$ with the corresponding p-values $\alpha$.*

Table 1
*Proportion of rejections of the null hypothesis at the 5% significance level in 10,000 simulations*

| $a_1$ | $-0.2$ | $-0.1$ | $0$ | $0.01$ |
|---|---|---|---|---|
| $\beta = 0$ | 0.014 | 0.026 | 0.049 | 0.052 |
| $\beta = 0.3, \sigma = 0.2$ | 0.066 | 0.115 | 0.215 | 0.224 |
| $\beta = 0.3, \sigma = 0.1$ | 0.188 | 0.301 | 0.439 | 0.451 |

Figure 3 implies a p-value of less than 1% for the existence of at least two modes, and a p-value of 7.33% for the existence of at least three modes.

Our second example concerns the detection of an increase in a failure rate. Gijbels and Heckman (2004) compare a global test and four versions of a localized test in a simulation study. A sample of size $m = 50$ is drawn from a distribution whose hazard rate $h(t)$ is modeled via $\log h(t) = a_1 \log t + \beta(2\pi\sigma^2)^{-1/2} \exp\{-(t-\mu)^2/(2\sigma^2)\}$. Table 1 shows the power of our procedure from Section 5 for the choices of parameters $a_1, \beta, \sigma$ used by Gijbels and Heckman (2004). The cases with $\beta = 0, a_1 \leq 0$ pertain to the null hypothesis of a nonincreasing failure rate, whereas $\beta = 0, a_1 = 0.01$ implies an increasing failure rate. The other eight cases result in a failure rate with a local increase. The power of the test introduced in Section 5 exceeds those of the five tests examined by Gijbels and Heckman (2004) in four of the nine cases that involve an increase in the failure rate.

## 7. Proofs.

7.1. *Proofs of Propositions* 2.1, 4.1 *and* 5.1. The proofs rely on the following elementary inequality, which we state without proof.

LEMMA 7.1. *Let $G_o$ and $G$ be distribution functions on an interval $(a,b)$ with densities $g_o$ and $g$, respectively. Suppose that $g - g_o \leq 0$ on $(a,c)$ and $g - g_o \geq 0$ on $(c,b)$, where $a < c < b$. Then, $G^{-1} \geq G_o^{-1}$.*

Note that the conditions in Lemma 7.1 are satisfied if, for instance, $g_o$ and $g$ are differentiable with derivatives satisfying $g' \geq g'_o$.

PROOF OF PROPOSITION 2.1. It is well known that $U_{(1)}, \ldots, U_{(n)}$ are distributed as the order statistics of $n$ independent random variables having uniform distribution on $[0,1]$. Suppose that $f$, and thus $f_{jk}$ is nondecreasing on $\mathcal{I}_{jk}$, where $k - j > 1$. The assumptions of Lemma 7.1 are then satisfied, with $g = f_{jk}$ and $g_o(x) := 1\{x \in \mathcal{I}_{jk}\}/|\mathcal{I}_{jk}|$. This implies that for $j < i < k$,

$$X_{(i)} = G^{-1}(U_{(i;j,k)}) \geq G_o^{-1}(U_{(i;j,k)}) = X_{(j)} + (X_{(k)} - X_{(j)})U_{(i;j,k)},$$

whence $T_{jk}(\mathbf{X}) \geq T_{jk}(\mathbf{U})$. In the case where $f$ is nonincreasing on $\mathcal{I}_{jk}$, the reverse inequality $T_{jk}(\mathbf{X}) \leq T_{jk}(\mathbf{U})$ follows from Lemma 7.1 with $g(x) = 1\{x \in \mathcal{I}_{jk}\}/|\mathcal{I}_{jk}|$ and $g_o := f_{jk}$. □

PROOF OF PROPOSITION 4.1. Again, we apply Lemma 7.1, this time with the densities

$$g(u) := |\mathcal{I}_{jk}| f_{jk}(X_{(j)} + |\mathcal{I}_{jk}|u)$$

and $g_o := g_S$ on $(0,1)$. Note that

$$\inf_{(0,1)} g' = |\mathcal{I}_{jk}|^2 \inf_{\mathcal{I}_{jk}} f'_{jk} = S \equiv g'_S.$$

It thus follows from Lemma 7.1 that

$$T_{jk}(\mathbf{X}) = \sum_{i=j+1}^{k-1} \beta(G^{-1}(U_{(i;j,k)})) \geq \sum_{i=j+1}^{k-1} \beta(G_S^{-1}(U_{(i;j,k)})).$$

As for the moments of $\beta(G_c^{-1}(U))$, first note that, generally,

$$\mathbb{E}h(\beta(G_c^{-1}(U))) = \int_0^1 h(\beta(u))(1 + c(u - 1/2))\, du$$

$$= \frac{1}{2}\int_{-1}^1 h(v)\left(1 + \frac{c}{2}v\right) dv$$

for $h:[-1,1] \to \mathbb{R}$. Letting $h(v) := v^j$ with $j = 1, 2$ shows that the first and second moments of $\beta(G_c^{-1}(U))$ are given by $c/6$ and $1/3$, respectively. Moreover, letting $h(v) := \exp(tv)$ yields

$$M_c(t) := \log \mathbb{E} \exp(t\beta(G_c^{-1}(U))) - ct/6 = \log(A(t) + cB(t)) - ct/6,$$

where

$$A(t) := \frac{1}{2} \int_{-1}^{1} e^{tv} \, dv = \sinh(t)/t = \sum_{k=0}^{\infty} \frac{t^{2k}}{(2k+1)!},$$

$$B(t) := \frac{1}{4} \int_{-1}^{1} e^{tv} v \, dv = (\cosh(t)/t - \sinh(t)/t^2)/2$$

$$= \frac{t}{6} \sum_{k=0}^{\infty} \frac{3}{2k+3} \frac{t^{2k}}{(2k+1)!}.$$

We have to show that $M_c(t) \leq t^2/6$ for any $t \neq 0$. To this end, note that $\partial M_c(t)/\partial c$ equals $B(t)/(A(t) + cB(t)) - t/6$ and $\partial^2 M_c(t)/\partial c^2 < 0$. Thus, $M_c(t)$ is strictly concave in $c \in \{c : A(t) + cB(t) > 0\}$. The equation $\partial M_c(t)/\partial c = 0$ is equivalent to $A(t) + cB(t)$ being equal to $6B(t)/t > 0$ and this means that $ct/6 = 1 - tA(t)/(6B(t))$. Hence, elementary manipulations of the series expansions yield

$$M_c(t) \leq \log\left(\frac{6B(t)}{t}\right) + \frac{tA(t)}{6B(t)} - 1$$

$$= \log\left(\sum_{k=0}^{\infty} \frac{3}{2k+3} \frac{t^{2k}}{(2k+1)!}\right)$$

$$+ \frac{t^2}{15} \sum_{k=0}^{\infty} \frac{5 \cdot 3}{(2k+5)(2k+3)} \frac{t^{2k}}{(2k+1)!} \bigg/ \sum_{k=0}^{\infty} \frac{3}{2k+3} \frac{t^{2k}}{(2k+1)!}$$

$$\leq \log\left(\sum_{k=0}^{\infty} \frac{(t^2/10)^k}{k!}\right) + \frac{t^2}{15}$$

$$= \frac{t^2}{6}. \qquad \square$$

PROOF OF PROPOSITION 5.1.  By construction, the vector $(X'_{(i)} - X'_{(0)})_{i=1}^{n+1}$ is distributed as the vector of order statistics of $n+1$ independent random variables with standard exponential distribution. Well-known facts imply that the variables $D'_i$ are independent with standard exponential distribution. Hence, $(W'_1, \ldots, W'_n) =_{\mathcal{L}} (U_{(1)}, \ldots, U_{(n)})$, while $W'_0 = 0$ and $W'_{n+1} = 1$.

We now assume that the failure rate is nondecreasing on $\mathcal{I}_{jk}$; the nonincreasing case is treated analogously. The function $G(x) := -\log(1 - F(x))$ is then convex on $\mathcal{I}_{jk}$. Hence $\alpha_s := D'_s/D_s$ is nondecreasing in $s \in \{j+1, \ldots, k\}$. Consequently, for $j < i < k$,

$$\begin{aligned}W_{(i;j,k)} - W'_{(i;j,k)} &= \frac{\sum_{s=j+1}^{i} D_s}{\sum_{s=j+1}^{k} D_s} - \frac{\sum_{s=j+1}^{i} \alpha_s D_s}{\sum_{s=j+1}^{k} \alpha_s D_s} \\ &= \frac{\sum_{s=j+1}^{i} \sum_{t=i+1}^{k} (\alpha_t - \alpha_s) D_s D_t}{\sum_{s=j+1}^{k} D_s \sum_{t=j+1}^{k} \alpha_t D_t} \\ &\geq 0.\end{aligned}$$

Hence, $T_{jk}(\mathbf{W}) \geq T_{jk}(\mathbf{W}')$. $\square$

7.2. *An auxiliary result concerning stochastic processes.* Our proof of Theorem 3.1 builds on a new version of Theorem 6.1 of Dümbgen and Spokoiny (2001). An important difference is that the original requirement of sub-Gaussian increments is relaxed to subexponential increments. The new version itself is just a corollary to more general results concerning stochastic processes obtained by Dümbgen and Walther (2006). We consider a stochastic process $Z = (Z(t))_{t \in \mathcal{T}}$ with continuous sample paths on a totally bounded metric space $(\mathcal{T}, \rho)$, where $\rho \leq 1$. "Totally bounded" means that for any $u > 0$, the capacity number

$$D(u, \mathcal{T}, \rho) := \sup\{\#\mathcal{T}_o : \mathcal{T}_o \subset \mathcal{T} \text{ such that } \rho(s,t) > v \text{ for different } s, t \in \mathcal{T}_o\}$$

is finite. In addition, we consider a function $\sigma : \mathcal{T} \to (0, 1]$, where $\sigma(t)$ may be viewed as a measure of spread of the distribution of $Z(t)$. We assume that

(7.1) $$|\sigma(s) - \sigma(t)| \leq \rho(s,t) \quad \text{for all } s, t \in \mathcal{T}$$

and that $\{t \in \mathcal{T} : \sigma(t) \geq \delta\}$ is compact for any $\delta \in (0, 1]$.

THEOREM 7.1. *Suppose that the following three conditions are satisfied:*

(i) *there exist constants $A, B, V > 0$ such that for arbitrary $u, \delta \in (0, 1]$,*

$$D(u\delta, \{t \in \mathcal{T} : \sigma(t) \leq \delta\}, \rho) \leq A u^{-B} \delta^{-V};$$

(ii) *there exists a constant $K \geq 1$ such that for arbitrary $s, t \in \mathcal{T}$ and $\eta \geq 0$,*

$$\mathbb{P}(|Z(s) - Z(t)| \geq K\rho(s,t)\eta) \leq K \exp(-\eta);$$

(iii) *for arbitrary $t \in \mathcal{T}$ and $\eta \geq 0$,*

$$\mathbb{P}(|Z(t)| \geq \sigma(t)\eta) \leq 2\exp(-\eta^2/2).$$

*Then,*

$$\mathbb{P}\left(\sup_{s,t\in\mathcal{T}}\frac{|Z(s)-Z(t)|}{\rho(s,t)\log(e/\rho(s,t))}\geq\eta\right)\leq p_1(\eta),$$

$$\mathbb{P}\left(\sup_{t\in\mathcal{T}}\frac{|Z(t)|/\sigma(t)-\sqrt{2V\log(1/\sigma(t))}}{D(\sigma(t))}\geq\eta\right)\leq p_2(\eta),$$

*with* $D(\delta) := \log(e/\delta)^{-1/2}\log(e\log(e/\delta))$, *where* $p_1 = p_1(\cdot|A,B,K)$ *and* $p_2 = p_2(\cdot|A,B,V,K)$ *are universal functions such that* $\lim_{\eta\to\infty} p_j(\eta) = 0$.

7.3. *Proof of Theorem* 3.1. In what follows, we describe a proof, but omit some technical arguments and details; for a complete account, we refer to Dümbgen and Walther (2006).

We embed our test statistics $T_{jk}$ into a stochastic process $Z_n$ on

$$\mathcal{T}_n := \{(\tau_{jn}, \tau_{kn}) : 0 \leq j < k \leq n+1\},$$

where $\tau_{in} := i/(n+1)$, equipped with the distance

$$\rho((u,v),(u',v')) := (|u-u'|+|v-v'|)^{1/2}$$

on $\mathcal{T} := \{(u,v): 0 \leq u < v \leq 1\}$. Namely, let

$$Z_n(\tau_{jn}, \tau_{kn}) := 3^{1/2}(n+1)^{-1/2} T_{jk}(\mathbf{U}).$$

Moreover, for $(u,v) \in \mathcal{T} \setminus \mathcal{T}_n$, let

$$Z_n(u,v) := Z_n(\tau_n(u), \tau_n(v)) \qquad \text{with } \tau_n(c) := \frac{\lfloor (n+1)c \rfloor}{n+1}.$$

Note that

$$\mathbb{E}(Z_n(u,v)) = 0 \quad \text{and} \quad \operatorname{Var}(Z_n(u,v)) \leq \sigma(u,v)^2,$$

where $\sigma(u,v) := (v-u)^{1/2}$. Elementary calculations show that these functions $\rho$ and $\sigma$ satisfy (7.1). Later, we shall prove the two following results concerning these processes $Z_n$ and the limiting process $Z$ defined in Theorem 3.1.

LEMMA 7.2. *The processes $Z$ on $\mathcal{T}$ and $Z_n$ on $\mathcal{T}_n$ ($n \in \mathbb{N}$) satisfy conditions* (i)–(iii) *of Theorem 7.1 with $A = 12$, $B = 4$, $V = 2$ and some universal constant $K$.*

LEMMA 7.3. *For any finite subset $\mathcal{T}_o$ of $\mathcal{T}$, the random variable $(Z_n(t))_{t\in\mathcal{T}_o}$ converges in distribution to $(Z(t))_{t\in\mathcal{T}_o}$.*

With arguments similar to those in Dümbgen (2002), one can deduce from Theorem 7.1 and Lemmas 7.2–7.3 that the preliminary test statistic

$$\tilde{T}_n := \max_{0 \le j < k \le n+1} \left( 3^{1/2}(k-j)^{-1/2} T_{jk}(\mathbf{U}) - \Gamma\left(\frac{k-j}{n+1}\right) \right)$$

$$= \max_{t \in \mathcal{T}_n} \left( \frac{|Z_n(t)|}{\sigma(t)} - \Gamma(\sigma(t)^2) \right)$$

[with $T_{jk}(\mathbf{U}) := 0$ if $k - j = 1$] converges in distribution to $T(W)$. Moreover,

$$T_n(\mathbf{U}) = \max_{t \in \mathcal{T}_n} \left( \frac{|Z_n(t)|}{\sigma_n(t)} - \Gamma(\sigma(t)^2) \right)$$

with

$$\sigma_n(t) := (\sigma(t)^2 - (n+1)^{-1})^{1/2},$$

where we use the convention that $0/0 := 0$. By means of the inequality $|Z_n(t)| \le (n+1)^{1/2} \sigma_n(t)^2$ and elementary considerations, one can show that $T_n(\mathbf{U}) = \tilde{T}_n + o_p(1)$, whence $T_n(\mathbf{U}) \to_{\mathcal{L}} T(W)$.

PROOF OF LEMMA 7.2. A proof of condition (i) is given by Dümbgen and Spokoiny (2001) (proof of Theorem 2.1) in a slightly different setting.

Next, we verify condition (ii). In order to bound the increment $Z_n(s) - Z_n(t)$ in terms of $\rho(s,t)$, we first consider the special case of $s = (0,1)$ and $t = (\tau, 1)$, where $\tau = \tau_{kn}$ for some $k \in \{1, \ldots, n\}$. Note that

$$\sum_{i=1}^{n} (2U_{(i)} - 1) = \sum_{i=1}^{k-1} (2U_{(i)} - 1) + 2U_{(k)} - 1 + \sum_{i=k+1}^{n} (2U_{(i)} - 1),$$

$$\sum_{i=1}^{k-1} (2U_{(i)} - 1) = \sum_{i=1}^{k-1} \left( 2\frac{U_{(i)}}{U_{(k)}} - 1 \right) U_{(k)} + (k-1)U_{(k)},$$

$$\sum_{i=k+1}^{n} (2U_{(i)} - 1) = \sum_{i=k+1}^{n} (2(U_{(i)} - U_{(k)}) - 1) + 2(n-k)U_{(k)}$$

$$= \sum_{i=k+1}^{n} \left( 2\frac{U_{(i)} - U_{(k)}}{1 - U_{(k)}} - 1 \right) (1 - U_{(k)}) + (n-k)U_{(k)},$$

whence

$$Z_n(0,1) = Z_n(0,\tau) U_{(k)} + Z_n(\tau, 1)(1 - U_{(k)}) + 3^{1/2}(n+1)^{1/2}(U_{(k)} - \tau).$$

Consequently,

$$Z_n(0,1) - Z_n(\tau, 1)$$
$$= (Z_n(0,\tau) - Z_n(\tau,1)) U_{(k)} + 3^{1/2}(n+1)^{1/2}(U_{(k)} - \tau)$$

$$= 3^{1/2}(n+1)^{-1/2}\left(\sum_{i=1}^{k-1}\beta\left(\frac{U_{(i)}}{U_{(k)}}\right) - \sum_{i=k+1}^{n}\beta\left(\frac{U_{(i)} - U_{(k)}}{1 - U_{(k)}}\right)\right)U_{(k)}$$

$$+ 3^{1/2}(n+1)^{1/2}(U_{(k)} - \tau)$$

$$=_{\mathcal{L}} 3^{1/2}(n+1)^{-1/2}\sum_{i=1}^{n-1}\beta(U_i')U_{(k)} + 3^{1/2}(n+1)^{1/2}(U_{(k)} - \tau),$$

where $U_1, \ldots, U_n, U_1', \ldots, U_{n-1}'$ are independent and identically distributed. Note that $U_{(k)}$ has a beta distribution with parameters $k$ and $n+1-k$. This entails that

$$\mathbb{P}\{\pm(U_{(k)} - \tau) \geq c\} \leq \exp(-(n+1)\Psi(\tau \pm c, \tau)) \qquad \text{for all } c \geq 0,$$

where $\Psi(x, \tau) := \tau \log(\tau/x) + (1-\tau)\log((1-\tau)/(1-x))$ if $x \in (0,1)$, and $\Psi(x, \tau) := \infty$ otherwise; see Proposition 2.1 of Dümbgen (1998). Elementary calculations show that $\Psi(\tau \pm c, \tau)$ is not smaller than $c^2/(2\tau(1-\tau) + 2c)$, whence

(7.2)
$$\mathbb{P}\{\pm(U_{(k)} - \tau) \geq c\} \leq \exp\left(-\frac{(n+1)c^2}{2\tau(1-\tau) + 2c}\right)$$

for all $c \geq 0$. Consequently, for any $r \geq 0$,

$$\mathbb{P}\{|3^{1/2}(n+1)^{1/2}(U_{(k)} - \tau)| \geq r\rho((0,1),(\tau,1))\}$$

$$= \mathbb{P}\{|3^{1/2}(n+1)^{1/2}(U_{(k)} - \tau)| \geq r\tau^{1/2}\}$$

$$= \mathbb{P}\left\{|U_{(k)} - \tau| \geq \frac{r\tau^{1/2}}{3^{1/2}(n+1)^{1/2}}\right\}$$

(7.3)
$$\leq 2\exp\left(-\frac{r^2\tau}{6\tau(1-\tau) + 12^{1/2}r(n+1)^{-1/2}\tau^{1/2}}\right)$$

$$\leq 2\exp\left(-\frac{r^2}{6 + 12^{1/2}r((n+1)\tau)^{-1/2}}\right)$$

$$\leq 2\exp\left(-\frac{r^2}{6 + 4r}\right)$$

$$\leq 4\exp(-r/4).$$

Here, we used the fact that $(n+1)\tau \geq 1$. Moreover, for any $r \geq 1$,

$$\mathbb{P}\left\{\left|3^{1/2}(n+1)^{-1/2}\sum_{i=1}^{n-1}\beta(U_i')U_{(k)}\right| \geq r\tau^{1/2}\right\}$$

$$\leq \mathbb{P}\left\{\left|(3/n)^{1/2}\sum_{i=1}^{n-1}\beta(U_i')\right| \geq r^{1/2}\right\} + \mathbb{P}\{U_{(k)} \geq r^{1/2}\tau^{1/2}\}$$

$$\leq 2\exp(-r/2) + \mathbb{P}\{U_{(k)} - \tau \geq r^{1/2}\tau^{1/2} - \tau\}$$

$$\leq 2\exp(-r/2) + \exp\left(-\frac{(n+1)(r^{1/2}-1)^2\tau}{2\tau(1-\tau) + 2(r^{1/2}-1)\tau^{1/2}}\right)$$

$$\leq 2\exp(-r/2) + \exp\left(-\frac{(n+1)(r^{1/2}-1)^2\tau^{1/2}}{2 + 2(r^{1/2}-1)}\right)$$

$$\leq 2\exp(-r/2) + \exp\left(-\frac{(n+1)^{1/2}(r^{1/2}-1)^2}{2r^{1/2}}\right).$$

Note that the probability in question is zero if $r$ is greater than $3^{1/2}(n+1)^{-1/2}(n-1)\tau^{-1/2}$ and the latter number is smaller than $3^{1/2}n$. Thus, suppose that $r \leq 3^{1/2}n$. Then,

$$\frac{(n+1)^{1/2}(r^{1/2}-1)^2}{2r^{1/2}} \geq \frac{(3^{-1/2}r+1)^{1/2}(r^{1/2}-1)^2}{2r^{1/2}} \geq 3^{-1}(r^{1/2}-1)^2.$$

Consequently, for all $r \geq 0$ and some positive constant $C_1$,

$$(7.4) \quad \mathbb{P}\left\{\left|3^{1/2}(n+1)^{-1/2}\sum_{i=1}^{n-1}\beta(U_i')U_{(k)}\right| \geq r\tau^{1/2}\right\} \leq C_1\exp(-r/C_1).$$

Combining (7.4) and (7.4) yields

$$(7.5) \quad \mathbb{P}\{|Z_n(0,1) - Z_n(\tau_{kn}, 1)| \geq r\rho((0,1),(\tau_{kn},1))\} \leq C_2\exp(-r/C_2)$$

for some positive constant $C_2$. Symmetry considerations show that the same bound applies to $s = (0,1)$ and $t = (0,\tau)$, that is,

$$(7.6) \quad \mathbb{P}\{|Z_n(0,1) - Z_n(0,\tau)| \geq r\rho((\tau_{kn},1),(0,1))\} \leq C_2\exp(-r/C_2).$$

In order to treat the general case, note that the processes $Z_n$ rescale as follows. For $0 \leq J < K \leq n+1$,

$$(Z_n(\tau_{J+j,n}, \tau_{J+k,n}))_{0 \leq j < k \leq K-J}$$
$$=_{\mathcal{L}} \sigma(\tau_{Jn}, \tau_{Kn})(Z_{K-J}(\tau_{j,K-J}, \tau_{k,K-J}))_{0 \leq j < k \leq K-J},$$

while for $0 \leq j < k \leq K - J$ and $0 \leq j' < k' \leq K - J$,

$$\rho((\tau_{J+j,n}, \tau_{J+k,n}), (\tau_{J+j',n}, \tau_{J+k',n}))$$
$$= \sigma(\tau_{Jn}, \tau_{Kn})\rho((\tau_{j,K-J}, \tau_{k,K-J}), (\tau_{j',K-J}, \tau_{k',K-J})).$$

With this rescaling, one can easily verify condition (ii) with $K = 2C_2$.

Finally, according to Proposition 4.1, $\mathbb{E}\exp(r\beta(U_i)) \leq \exp(r^2/6)$ for all $r \in \mathbb{R}$, whence

$$\mathbb{E}\exp(r\sigma(t)^{-1}Z_n(t)) \leq \exp(r^2/2) \qquad \text{for } r \in \mathbb{R}, t \in \mathcal{T}_n.$$

A standard argument involving Markov's inequality then yields condition (iii). □

PROOF OF LEMMA 7.3. Recall the representation $U_{(i)} - U_{(i-1)} = E_i/S_n$ with independent, standard exponential variables $E_i$ and $S_n = \sum_{j=1}^{n+1} E_j$. Starting from (2.1), one can write

$$Z_n(\tau_{jn}, \tau_{kn}) = -3^{1/2} \frac{k-j}{(U_{(k)} - U_{(j)})S_n}(n+1)^{-1/2} \sum_{i=1}^{n+1} \beta\left(\frac{i-j-1/2}{k-j}\right) E_i$$

$$= \frac{\tau_{kn} - \tau_{jn}}{U_{(k)} - U_{(j)}} \frac{n+1}{S_n} \times \tilde{Z}_n(\tau_{jn}, \tau_{kn}),$$

where

$$\tilde{Z}_n(\tau_{jn}, \tau_{kn}) := 3^{1/2}(n+1)^{-1/2} \sum_{i=1}^{n+1} \beta\left(\frac{\tau_{in} - \tau_{jn} - \delta_n}{\tau_{kn} - \tau_{jn}}\right)(1 - E_i)$$

and $\delta_n := (2(n+1))^{-1}$. The centering of the variables $E_i$ is possible because the sum of the coefficients $\beta((i-j-1/2)/(k-j))$, $j < i \le k$, is zero. Since $S_n/(n+1) \to_p 1$ and $\max_{1 \le i \le n} |U_{(i)} - \tau_{in}| \to_p 0$, it suffices to consider the stochastic process $\tilde{Z}_n$ in place of $Z_n$. But the assertion then follows from the multivariate version of Lindeberg's central limit theorem and elementary covariance calculations. □

7.4. *Proofs for Section 4.* We first prove the lower bounds comprising Theorem 4.1 (b)–(c). The following lemma is a surrogate for Lemma 6.2 of Dümbgen and Spokoiny (2001) in order to treat likelihood ratios and i.i.d. data.

LEMMA 7.4. *Let $X_1, X_2, \ldots, X_n$ be i.i.d. with distribution $P$ on some measurable space $\mathcal{X}$. Let $f_1, \ldots, f_m$ be probability densities with respect to $P$ such that the sets $B_j := \{f_j \ne 1\}$ are pairwise disjoint and define $L_j := \prod_{i=1}^n f_j(X_i)$. Then,*

$$\mathbb{E}\left|m^{-1} \sum_{j=1}^m L_j - 1\right| \to 0,$$

*provided that $m \to \infty$, $\Delta_\infty \le C(\log m)^{-1/2}$ for some fixed constant $C$ and*

$$\sqrt{\log m}\left(1 - \frac{n\Delta_2^2}{2\log m}\right) \to \infty,$$

*where $\Delta_\infty := \max_j \sup_x |f_j(x) - 1|$ and $\Delta_2 := \max_j(\int (f_j - 1)^2 \, dP)^{1/2}$.*

PROOF OF LEMMA 7.4. The likelihood ratio statistics $L_j$ are not stochastically independent, but conditional on $\boldsymbol{\nu} = (\nu_j)_{j=1}^m$ with $\nu_j := \#\{i : X_i \in B_j\}$, they are. Furthermore, $\mathbb{E}(L_j) = 1 = \mathbb{E}(L_j|\boldsymbol{\nu})$. Thus, a standard truncation argument shows that for any $\epsilon > 0$ and $0 < \gamma \le 1$,

$$\mathbb{E}\left(\left|m^{-1}\sum_j L_j - 1\right|\bigg|\boldsymbol{\nu}\right)$$

$$\le m^{-1}\operatorname{Var}\left(\sum_j 1\{L_j \le \epsilon m\}L_j \bigg|\boldsymbol{\nu}\right)^{1/2} + 2m^{-1}\sum_j \mathbb{E}(1\{L_j > \epsilon m\}L_j|\boldsymbol{\nu})$$

$$\le m^{-1}\left(\sum_j \mathbb{E}(1\{L_j \le \epsilon m\}L_j^2|\boldsymbol{\nu})\right)^{1/2} + 2m^{-1}\sum_j \mathbb{E}(1\{L_j > \epsilon m\}L_j|\boldsymbol{\nu})$$

$$\le m^{-1}\left(\sum_j \mathbb{E}(\epsilon m L_j|\boldsymbol{\nu})\right)^{1/2} + 2\epsilon^{-\gamma}m^{-(1+\gamma)}\sum_j \mathbb{E}(L_j^{1+\gamma}|\boldsymbol{\nu})$$

$$= \epsilon^{1/2} + 2\epsilon^{-\gamma}m^{-(1+\gamma)}\sum_j \mathbb{E}(L_j^{1+\gamma}|\boldsymbol{\nu}).$$

It thus suffices to show that

$$\inf_{\gamma \in (0,1]} \max_j m^{-\gamma}\mathbb{E}(L_j^{1+\gamma}) \to 0$$

under the stated conditions on $m$, $\Delta_\infty$ and $\Delta_2$. Note that $\mathbb{E}(L_j^{1+\gamma})$ is equal to $\mathbb{E}(f_j(X_1)^{1+\gamma})^n$ and that elementary calculus gives

$$(1+y)^{1+\gamma} \le 1 + (1+\gamma)y + \gamma(1+\gamma)y^2/2 + 3\gamma|y|^3 \quad \text{for } |y| \le 1.$$

Hence, $\mathbb{E}(f_j(X_1)^{1+\gamma}) \le 1 + \gamma(1+\gamma)\Delta_2^2/2 + 3\gamma\Delta_\infty\Delta_2^2$ and

$$\max_j m^{-\gamma}\mathbb{E}(L_j^{1+\gamma})$$

(7.7)
$$\le m^{-\gamma}(1 + \gamma(1+\gamma)\Delta_2^2/2 + 3\gamma\Delta_\infty\Delta_2^2)^n$$
$$\le \exp(-\gamma\log m + \gamma(1+\gamma)n\Delta_2^2/2 + 3\gamma\Delta_\infty n\Delta_2^2).$$

Suppose that $n\Delta_2^2 \le 2(1-b_m)\log m$, where $(0,1) \ni b_m \to 0$ and $b_m^2\log m \to \infty$ as $m \to \infty$. Then, the right-hand side of (7.7) does not exceed

$$\exp(-\gamma(1-(1+\gamma)(1-b_m))\log m + 6\gamma\Delta_\infty \log m)$$

$$\le \exp\left(-\frac{b_m^2\log m}{4(1-b_m)} + 3Cb_m(\log m)^{1/2}\right) \quad \text{if } \gamma = \frac{b_m}{2(1-b_m)}$$

$$\to 0 \quad \text{as } m \to \infty. \qquad \square$$

PROOF OF THEOREM 4.1(b). Let $\tilde{c}_n := c_n\sqrt{\log(e/\delta_n)/n}$, and set $f_0 := 1_{[0,1)}$ and

$$f_{nj}(x) := f_0(x) + 1\{x \in I_{nj}\}\tilde{c}_n\delta_n^{-3/2}(x - (j-1/2)\delta_n)$$

for $j = 1, \ldots, m_n := \lfloor 1/\delta_n \rfloor$ and $I_{nj} := [(j-1)\delta_n, j\delta_n)$. Each $f_{nj}$ is a probability density with respect to the uniform distribution on $[0,1)$ such that the corresponding distribution $F_{nj}$ satisfies $F_{nj}(I_{nj}) = \delta_n$ and $\inf_{I_{nj}} f'_{nj} \cdot |I_{nj}|^2/\sqrt{F_{nj}(I_{nj})} = \tilde{c}_n$, that is, $f_{nj} \in \mathcal{F}(\delta_n, \tilde{c}_n)$. Thus, for any test $\phi_n(\mathbf{X})$ with $\mathbb{E}_{f_0}\phi_n(\mathbf{X}) \leq \alpha + o(1)$,

$$\inf_{f \in \mathcal{F}(\delta_n, \tilde{c}_n)} \mathbb{E}_f \phi_n(\mathbf{X}) - \alpha \leq m_n^{-1} \sum_{j=1}^{m_n} \mathbb{E}_{f_{nj}} \phi_n(\mathbf{X}) - \alpha$$

$$= \mathbb{E}_{f_0}\left(\left(m_n^{-1}\sum_{j=1}^{m_n} L_{nj} - 1\right)\phi_n(\mathbf{X})\right) + o(1)$$

$$\leq \mathbb{E}_{f_0}\left|m_n^{-1}\sum_{j=1}^{m_n} L_{nj} - 1\right| + o(1),$$

where $L_{nj} := \prod_{i=1}^n f_{nj}(X_i)$. The latter expectation tends to zero by Lemma 7.4. For $\Delta_2^2 = \tilde{c}_n^2/12$, and $\Delta_\infty = \tilde{c}_n\delta_n^{-1/2}/2$ is less than $\sqrt{6\log(e/\delta_n)/(n\delta_n)} = O(\log(n)^{-1/2}) = O(\log(m_n)^{-1/2})$ because $n\delta_n \geq \log(n)^2$, hence $m_n = \delta_n^{-1} + O(1) = o(n)$. Finally,

$$\sqrt{\log m_n}\left(1 - \frac{n\Delta_2^2}{2\log m_n}\right) = \frac{24\log m_n - c_n^2 \log(e/\delta_n)}{24\sqrt{\log m_n}}$$

$$\geq \sqrt{24}(\sqrt{24} - c_n)\sqrt{\log(e/\delta_n)}(1 + o(1)) + o(1)$$

tends to infinity by assumptions on $\delta_n$ and $c_n$. □

PROOF OF THEOREM 4.1(c). We may assume without loss of generality that the left endpoint of $I_n$ is 0. We now define probability densities $f_n$ and $g_n$ via

$$f_n(x) := \frac{\delta_n}{|I_n|}1\{x \in [0, |I_n|/\delta_n]\},$$

$$g_n(x) := f_n(x) + \frac{\sqrt{\delta_n}b_n}{\sqrt{n}|I_n|^2}(x - |I_n|/2)1\{x \in I_n\}.$$

Note that $g_n \geq 0$ because $b_n \leq 2\sqrt{n\delta_n}$. Furthermore, $f_n$ is nonincreasing on $I_n$, while $g_n$ belongs to $\mathcal{F}(I_n, \delta_n, b_n/\sqrt{n})$.

We now apply LeCam's notion of contiguity [cf. LeCam and Yang (1990), Chapter 3]: If a test $\phi_n(\mathbf{X})$ satisfies $\mathbb{E}_{f_n}\phi_n(\mathbf{X}) \leq \alpha$, then $\limsup \mathbb{E}_{g_n}\phi_n(\mathbf{X}) < 1$, provided that

$$(7.8) \qquad \mathcal{L}_{f_n}\left(\sum_{i=1}^n \log(g_n/f_n)(X_i)\right) \to_w Q$$

for some probability measure $Q$ on the real line such that $\int e^x Q(dx) = 1$.

Note that $\mathcal{L}_{f_n}(\sum_{i=1}^n \log(g_n/f_n)(X_i))$ equals the distribution of $\sum_{i=1}^{N_n} \log(1 + c_n V_i)$ with $c_n := b_n/(2\sqrt{n\delta_n}) \in [0,1]$ and independent random variables $N_n$, $V_1, V_2, V_3, \ldots$ such that $N_n \sim \mathrm{Bin}(n, \delta_n)$ and $V_i \sim \mathrm{Unif}[-1,1]$.

First, suppose that $n\delta_n \not\to \infty$. By extracting a subsequence, if necessary, we may assume that $n\delta_n \to \lambda \in [0, \infty)$ and $c_n \to c \in [0,1]$. (7.8) then holds for the distribution $Q := \sum_{k=0}^\infty p_\lambda(k) \mathcal{L}(\sum_{i=1}^k \log(1 + cV_i))$ with the Poisson weights $p_\lambda(k) := e^{-\lambda}\lambda^k/k!$. But this measure $Q$ satisfies $\int e^x Q(dx) = 1$, whence $\limsup \mathbb{E}_{g_n}\phi_n(\mathbf{X}) < 1$. This contradiction shows that $n\delta_n \to \infty$.

Second, suppose that $n\delta_n \to \infty$, but $b_n \not\to \infty$. We assume without loss of generality that $b_n \to b \in [0, \infty)$. Lindeberg's central limit theorem and elementary calculations yield (7.8) with Gaussian distribution $Q = \mathcal{N}(-b^2/24, b^2/12)$. Again, the limit distribution satisfies $\int e^x Q(dx) = 1$. Hence, $b_n \to \infty$. □

Theorem 4.2 concerns our specific multiscale procedure. It will be derived from the following basic result.

LEMMA 7.5. *For a bounded open interval $I$ and $\delta \in (0,1]$ let $f$ be a density in $\mathcal{F}(I, \delta, D\sqrt{\log(e/\delta)/n})$ with $D \geq \sqrt{24}$. Then,*

$$n\delta \geq \tilde{D}\max(\log(e/\delta), K\log(en)),$$

*with $\tilde{D} := D^2/4$ and $K \geq 1 - (\log \tilde{D} + \log\log(en))/\log(en)$. Suppose that*

$$(7.9) \quad D \geq \frac{\sqrt{24}}{(1-\epsilon)^2\sqrt{1-\gamma-2/(n\delta)}}\left(1 + \frac{\kappa_n(\alpha)+\eta}{\Gamma(\delta)} + \frac{\gamma + 2/(n\delta)}{\Gamma(\delta)^2}\right)$$

*for certain numbers $\epsilon \in (0,1)$, $\gamma \in (0, 1/2]$ and $\eta > 0$. Then*

$\mathbb{P}(\mathcal{D}^+(\alpha)$ *contains no interval $J \subset I$)*

$$\leq \exp(-n\delta\gamma^2/2) + 2\exp(-D\sqrt{n\delta\log(e/\delta)}\epsilon^2/8) + \exp(-\eta^2/2).$$

PROOF. The inequalities $2\sqrt{\delta} \geq H(f, I) \geq D\sqrt{\log(e/\delta)/n}$ entail that $n\delta \geq \tilde{D}\log(e/\delta)$. Now, write $n\delta = \tilde{D}K\log(en)$ for some $K > 0$. In the case of $K \leq 1$,

$$\tilde{D}K\log(en) \geq \tilde{D}\log(e/\delta) = \tilde{D}(\log(en) - \log(\tilde{D}K\log(en)))$$

$$\geq \tilde{D}\log(en)\left(1 - \frac{\log\tilde{D} + \log\log(en)}{\log(en)}\right),$$

and dividing both sides by $\tilde{D}\log(en)$ yields the asserted lower bound for $K$.

The number $N := \#\{i : X_{(i)} \in I\}$ has distribution $\mathrm{Bin}(m,\delta)$ with $m \in \{n, n+1, n+2\}$. Consequently it follows from Chernov's exponential inequality for binomial distributions [cf. van der Vaart and Wellner (1996), A.6.1] that

$$\mathbb{P}(N \leq (1-\gamma)n\delta) \leq \exp(-n\delta\gamma^2/2).$$

Since $D \geq \sqrt{24}$ by assumption, we can conclude that $n\delta \geq D^2/4 > 6$, so $(1-\gamma)n\delta \geq 3$. In the case of $N \geq 3$, let $j := \min\{i : X_{(i)} \in I\}$ and $k := \max\{i : X_{(i)} \in I\}$, that is, $N = k - j + 1$. In order to bound the probability of $|\mathcal{I}_{jk}|/|I| < 1 - \epsilon$, we write $I = (a,b)$ and define $I_{(\ell)} := (a, a + \epsilon|I|/2]$, $I_{(r)} := [b - \epsilon|I|/2, b)$. Then,

$$nF(I_{(r)}) \geq nF(I_{(\ell)}) \geq n\inf_I f' \cdot |I_{(\ell)}|^2/2 = nH(f,I)\sqrt{\delta}\epsilon^2/8$$
$$\geq D\sqrt{n\delta\log(e/\delta)}\epsilon^2/8,$$

whence

$$\mathbb{P}(N \leq 1 \text{ or } |\mathcal{I}_{jk}|/|I| \leq 1 - \epsilon)$$
$$\leq \mathbb{P}(\text{no observations in } I_{(\ell)}) + \mathbb{P}(\text{no observations in } I_{(r)})$$
$$\leq 2\exp(-D\sqrt{n\delta\log(e/\delta)}\epsilon^2/8).$$

Hereafter, we will always assume that $N \geq (1-\gamma)n\delta$ and $|\mathcal{I}_{jk}|/|I| \geq 1-\epsilon$. By $\mathbb{P}^*(\cdot)$, we denote conditional probabilities given these two inequalities. The definition of $\mathcal{D}^+(\alpha)$ implies that $\mathbb{P}^*(\mathcal{D}^+(\alpha)$ contains no $J \subset I)$ is not greater than $\mathbb{P}^*(T_{jk}(\mathbf{X}) \leq c_{jk}(\alpha))$. On the other hand, it follows from Proposition 4.1 that

$$\mathbb{P}^*\left(T_{jk}(\mathbf{X}) \leq \frac{\tilde{C}(N-2)}{6} - \eta\sqrt{\frac{N-2}{3}}\right) \leq \exp(-\eta^2/2) \qquad \text{for any } \eta \geq 0,$$

where $\tilde{C} := H(f, \mathcal{I}_{jk})/\sqrt{F(\mathcal{I}_{jk})}$. It thus suffices to show that

$$\frac{\tilde{C}(N-2)}{6} - \eta\sqrt{\frac{N-2}{3}} \geq c_{jk}(\alpha).$$

By the definition of $c_{jk}(\alpha)$ this is equivalent to

$$\tilde{C}\sqrt{\frac{N-2}{12}} \geq \Gamma\left(\frac{N-1}{n+1}\right) + \kappa_n(\alpha) + \eta.$$

However, the left-hand side is not smaller than

$$\frac{(1-\epsilon)^2 H(f,I)}{\sqrt{\delta}}\sqrt{\frac{N-2}{12}} \ge \frac{(1-\epsilon)^2 H(f,I)\sqrt{(1-\gamma)n\delta - 2}}{\sqrt{12\delta}}$$

$$\ge D\frac{(1-\epsilon)^2\sqrt{1-\tilde\gamma}}{\sqrt{24}}\Gamma(\delta)$$

$$\ge \Gamma(\delta) + \kappa_n(\alpha) + \eta + \frac{\tilde\gamma}{\Gamma(\delta)}$$

with $\tilde\gamma := \gamma + 2/(n\delta)$, whereas $\Gamma((N-1)/(n+1)) \le \Gamma((N-2)/n)$ is not greater than

$$\Gamma(\delta(1-\tilde\gamma)) \le \Gamma(\delta) - \log(1-\tilde\gamma)/\Gamma(\delta) \le \Gamma(\delta) + \tilde\gamma/\Gamma(\delta). \qquad \square$$

PROOF OF THEOREM 4.2. First note that (4.2) and the first part of Lemma 7.5 entail that

$$n\delta_n \ge C^2/4 \ge 6 \quad \text{and} \quad n\delta_n \ge (C^2/4 + o(1))\log n.$$

In particular, $\#\mathcal{I}_n \le \delta_n^{-1} = o(n)$.

We apply Lemma 7.5 to $f = f_n$ and all intervals $I \in \mathcal{I}_n$. More precisely, we shall introduce suitable numbers $\gamma_n \in (0, 1/2]$, $\epsilon_n \in (0,1)$ and $\eta_{n,I} > 0$. According to Lemma 7.5, the probability that some $I \in \mathcal{I}_n$ does not cover an interval from $\mathcal{D}^+(\alpha)$ is bounded by

$$\begin{aligned}(7.10) \quad &\#\mathcal{I}_n(\exp(-n\delta_n\gamma_n^2/2) + 2\exp(-C\sqrt{n\delta_n\log(e/\delta_n)}\epsilon_n^2/8)) \\ &+ \sum_{I\in\mathcal{I}_n}\exp(-\eta_{n,I}^2/2),\end{aligned}$$

provided that

$$C\left(1 + \frac{\sqrt{2}b_n}{\Gamma(F_n(I))}\right) \ge \frac{\sqrt{24}}{(1-\epsilon_n)^2\sqrt{1-\tilde\gamma_n}}\left(1 + \frac{\kappa_n(\alpha) + \eta_{n,I}}{\Gamma(F_n(I))} + \frac{\tilde\gamma_n}{\Gamma(F_n(I))^2}\right)$$

for all $I \in \mathcal{I}_n$, where $\tilde\gamma_n := \gamma_n + 2/(n\delta_n) = O(1)$. Also note that $\kappa_n(\alpha) = O(1)$, by virtue of Theorem 3.1. Hence, the preceding requirement is met if, for every constant $A > 0$ and sufficiently large $n$,

$$(7.11) \quad C\left(1 + \frac{\sqrt{2}b_n}{\Gamma(F_n(I))}\right) \ge \frac{\sqrt{24}}{(1-\epsilon_n)^2\sqrt{1-\tilde\gamma_n}}\left(1 + \frac{A + \eta_{n,I}}{\Gamma(F_n(I))}\right)$$

for all $I \in \mathcal{I}_n$.

In setting (i), we use constants $\gamma_n = \gamma \in (0, 1/2]$ and $\epsilon_n = \epsilon \in (0,1)$ to be specified later and define

$$\eta_{n,I} := \sqrt{2\log(1/F_n(I)) + b_n} \le \Gamma(F_n(I)) + \sqrt{b_n}.$$

Since $\delta \log(e/\delta)$ is nondecreasing in $\delta \in (0,1]$, it follows from $n\delta_n \geq (C^2/4 + o(1)) \log n$ that

$$\sqrt{n\delta_n \log(e/\delta_n)} \geq (C/2 + o(1)) \log n.$$

Hence, the bound in (7.10) equals

$$o(1) \cdot (\exp(-(C^2\gamma^2/8 - 1 + o(1)) \log n)$$
$$+ \exp(-(C^2\epsilon^2/16 - 1 + o(1)) \log n)) + \sum_{I \in \mathcal{I}_n} F_n(I) \exp(-b_n/2)$$

and tends to zero, provided that $\gamma > \sqrt{8}/C$ and $\epsilon > 4/C$. Moreover, the right-hand side of (7.11) is not greater than

$$\frac{\sqrt{24}}{(1-\epsilon)^2 \sqrt{1-\gamma - o(1)}} \left(2 + \frac{A + \sqrt{b_n}}{\Gamma(F_n(I))}\right)$$
$$= \frac{2\sqrt{24} + o(1)}{(1-\epsilon^2)\sqrt{1-\gamma}} \left(1 + \frac{o(b_n)}{\Gamma(F_n(I))}\right).$$

Hence, the conclusion for setting (i) is correct if, say, $\epsilon = 4/(2\sqrt{24}) = \sqrt{1/6}$ and $\gamma = \sqrt{8}/(2\sqrt{24}) = \sqrt{1/12}$, while $C$ is strictly larger than

$$\frac{2\sqrt{24}}{(1-\epsilon)^2\sqrt{1-\gamma}} < 34.$$

In settings (ii)–(iii), we define

$$\gamma_n := (2\log(D\#\mathcal{I}_n)/(n\delta_n))^{1/2},$$
$$\epsilon_n := ((8/C)\log(D\#\mathcal{I}_n)/\sqrt{n\delta_n \log(e/\delta_n)})^{1/2},$$
$$\eta_{n,I} := \begin{cases} \sqrt{2\log(1/F_n(I)) + b_n}, & \text{in setting (ii)}, \\ b_n/D, & \text{in setting (iii)}, \end{cases}$$

for some (large) constant $D > 1$. The bound in (7.10) is then not greater than

$$3/D + \begin{cases} \exp(-b_n/2), & \text{in setting (ii)} \\ \exp(\log \#\mathcal{I}_n - b_n^2/(2D^2)), & \text{in setting (iii)} \end{cases} = 3/D + o(1).$$

It thus remains to verify (7.11).

Note that $\gamma_n \to 0$, by assumption. Moreover, since $\#\mathcal{I}_n \leq \delta_n^{-1}$, the term $\log(D\#\mathcal{I}_n)$ is not greater than $\log(D/\delta_n)^{1/2} \log(D\#\mathcal{I}_n)^{1/2}$, whence

$$\epsilon_n \leq \sqrt{8/C}(\log D)^{1/4}(\log(D\#\mathcal{I}_n)/(n\delta_n))^{1/4} \to 0.$$

Hence, in setting (ii), the right-hand side of (7.11) is not greater than

$$(2\sqrt{24} + o(1))\left(1 + \frac{o(b_n)}{\Gamma(F_n(I))}\right),$$

so (7.11) is satisfied for sufficiently large $n$, if $C > 2\sqrt{24}$. In setting (iii), the right-hand side of (7.11) is not greater than

$$\sqrt{24}\left(1 + O(\tilde{\gamma}_n + \epsilon_n) + (1 + o(1))\frac{A + b_n/D}{\Gamma(F_n(I))}\right).$$

By the first part of Lemma 7.5, $n\delta_n \geq (C^2/4)\log(e/\delta_n) \geq (C^2/8)\Gamma(F_n(I))^2$ for all $I \in \mathcal{I}_n$. Thus,

$$\tilde{\gamma}_n + \epsilon_n \leq \frac{O(\log(D\#\mathcal{I}_n)^{1/2})}{\Gamma(F_n(I))} = \frac{o(b_n)}{\Gamma(F_n(I))} \qquad \text{for all} I \in \mathcal{I}_n.$$

Consequently, (7.11) is satisfied if $C \geq \sqrt{24}$. $\square$

Institute of Mathematical Statistics
and Actuarial Science
University of Bern
Sidlerstrasse 5
CH-3012 Bern
Switzerland
E-mail: duembgen@stat.unibe.ch

Department of Statistics
Stanford University
390 Serra Mall
Stanford, California 94305
USA
E-mail: walther@stat.stanford.edu